\documentclass{article}

\usepackage[ansinew]{inputenc}
\usepackage[english]{babel}
\usepackage{graphicx}
\usepackage{amsfonts,amssymb,amsmath,amsthm}

\newtheorem {thm}{Theorem}

\newtheorem {defn}[thm]{Definition}

\newtheorem {mthm}[thm]{Main Theorem}
\newtheorem {lemma}[thm]{Lemma}

\newtheorem {cor}[thm]{Corollary}

\newcommand{\TR}{\operatorname{tr}}

\newcommand{\bbmatrix}[1]{\left[ \begin{array}{ccccccccccccccccccc} #1 \end{array} \right]}

\newcommand{\Rem}{\noindent\textbf{Remark: }}

\newcommand{\Proof}{\noindent\textbf{Proof: }}

\title{Classification of KdV vessels with constant parameters and
two dimensional outer space}
\author{Andrey Melnikov\\SunEdison, MO, USA}

\begin{document}

\maketitle
\abstract{In this article we classify vessels producing solutions of some completely integrable PDEs, presenting a \textit{unified} approach for them. The classification includes such important examples as Korteweg-de Vries (KdV) and evolutionary Non Linear Schr\" odingier (ENLS) equations. In fact, employing basic matrix algebra techniques it is shown that there are exactly two canonical forms of such vessels, so that each canonical form generalize either KdV or ENLS equations. Particularly, Dirac canonical systems, whose evolution was recently inserted into the vessel theory, are shown to be equivalent to the ENLS equation in the sense of vessels. This work is important as a first step to classification of completely integrable PDEs, which are solvable by the theory of vessels. We note that a recent paper of the author, published in Journal of Mathematical Physics, showed that initial value problem with analytic initial potential for the KdV equation has at least a "narrowing" in time solution. The presented classification, inherits this idea and a similar theorem can be easily proved for the presented PDEs. Finally, the the resuts of the work serve as a basis for the investigation of the following problems: 1. hierarchy of the generalized KdV, ENLS equations (by generalizing the vessel equations), 2. new completely integrable PDEs (by changing the dimension of the outer space), 3. addressing the question of integrability of a given arbitrary PDE (the future classification will create a list of solvable by vessels equations, which may eventually include many existing classes of PDEs).}

\section{Background}
In  series of recent papers of the author, it was shown that some completely integrable PDEs can be fully addressed using a theory of vessels. For example, the famous Korteweg-de Vries (KdV) equation \cite{bib:KdV, bib:GGKM}, which was first considered by the author in \cite{bib:GenVessel}, was further significantly improved in \cite{bib:KdVVessels} and finally culminated in a scattering theory of analytic parameters in  \cite{bib:AnalPotKdV}, and KdV hierarchy \cite{bib:KdVHierarchy}. Second example is a notable Evolutionary Non Linear Schr\" odinger (ENLS) equation \cite{bib:KatoNLS}, which was inserted into the setting of vessels in \cite{bib:ENLS}. Third example, developed by the author addressed the scattering theory and a corresponding completely integrable PDE for so called canonical (or Dirac) systems \cite{bib:CanSys}. Worth noticing that a different from these three examples type of vessels was used in \cite{bib:BoussVessels} to solve a shallow waters (Boussinesq) equation, originated in \cite{bib:Bous}.

The scheme for construction of solutions of the KdV, ENLS and canonical PDE equations is the same. It is a generalization of Zacharov-Shabath scheme written in a matricial form. This scheme is further elaborated into a more detailed notion of so called KdV vessels, which in a very special case solve the KdV equation, mentioned earlier. More precisely, we consider the following collection of bounded operators and spaces ($\mathcal K$ is a Krein space)
\[ \begin{array}{lllllllll}
\mathfrak{V}_{KdV} = \bbmatrix{C(x,t) & A_\zeta, \mathbb X(x,t), A & B(x,t)&  \sigma_1,\sigma_2,\gamma,\gamma_*(x,t) \\
\mathcal E & \mathcal{K} & \mathcal E & \Omega}, \\
\hspace{2cm} B, C^*:\mathbb C^2\rightarrow\mathcal K, \quad A,\mathbb X,A_\zeta:\mathcal K\rightarrow\mathcal K, \\
\hspace{3cm} \sigma_1, \sigma_2, \gamma,\gamma_*(x, t) \text{ - $2\times 2$ matrices}, \quad \Omega\subseteq\mathbb R^2\\
\left\{ \begin{array}{ll}
\dfrac{\partial}{\partial x} B \sigma_1 = - A B \sigma_2 - B \gamma), \\
\dfrac{\partial}{\partial t} B = i A  \dfrac{\partial}{\partial x} B,
\end{array} \right. \quad \quad 
\left\{ \begin{array}{ll}
\sigma_1 \dfrac{\partial}{\partial x} C   = - \sigma_2 C A_\zeta + \gamma C, \\ 
\dfrac{\partial}{\partial t} C  = - i \dfrac{\partial}{\partial x} C A_\zeta,
\end{array}\right. \\
\hspace{3cm} \left\{ \begin{array}{ll}
\dfrac{\partial}{\partial x} \mathbb X =  B \sigma_2 C, \\
 \dfrac{\partial}{\partial t} \mathbb X = i A B \sigma_2 C - i B \sigma_2 C A_\zeta + i B \gamma C, \\
\end{array}\right. \\
A \mathbb X + \mathbb X A_\zeta  + B \sigma_1 C  = 0, \\
\gamma_* =  \gamma + \sigma_2 C \mathbb X^{-1} B \sigma_1 - \sigma_1 C \mathbb X^{-1} B \sigma_2.
\end{array} \]
which is called a \textbf{(regular) KdV vessel} on $\Omega$, where $\mathbb X(x, t)$ is invertible. The author showed in
\cite{bib:GenVessel, bib:SLVessels, bib:CanSys, bib:ENLS, bib:AnalPotKdV} that in each particular case these equations are easily solvable. To see it, notice that the the system of equations for $B(x,t)$
\[
\left\{ \begin{array}{ll}
\dfrac{\partial}{\partial x} B \sigma_1 = - (A B \sigma_2 + B \gamma), \\
\dfrac{\partial}{\partial t} B = i A  \dfrac{\partial}{\partial x} B,
\end{array} \right.  \]
consists of a constant coefficients equation (upper), and an easily solvable wave equation (bottom). Similarly, the system of equations for $C(x,t)$ is solved. Finally, the operator $\mathbb X(x,t)$ is obtained from the operators  $B(x,t), C(x,t)$ by integration. The last two equations are also easily handled: in fact, the last equation serves as a defining equation for the matrix $\gamma_*(x,t)$ and the algebraic relation $A \mathbb X + \mathbb X A_\zeta  + B \sigma_1 C  = 0$ must be satisfied as an initial condition only (or more precisely, for fixed $x_0,t_0$).

To obtain a solution of the KdV equation \eqref{eq:KdV} on a set $\Omega$ we choose the outer space
$\mathcal E=\mathbb C^2$ and
\[ \sigma_1 = \bbmatrix{0&1\\1&0}, \quad \sigma_2=\bbmatrix{1&0\\0&0}, \quad \gamma=\bbmatrix{0&0\\0&i}.
\]
Then the function ($(x_0,t_0)\in\Omega$)
\[ q(x,t) = -2\dfrac{\partial^2}{\partial x^2} \ln \det(\mathbb X^{-1}(x_0,t_0)\mathbb X(x,t))
\]
satisfies \cite{bib:KdVVessels, bib:AnalPotKdV}
\begin{equation} \label{eq:KdV}
q_t = -\dfrac{3}{2} q q_x + \dfrac{1}{4} q_{xxx}.
\end{equation}
If we choose vessel parameters 
\[ \sigma_1 = I= \bbmatrix{1&0\\0&1}, \quad \sigma_1 = \dfrac{1}{2}\bbmatrix{1&0\\0&-1}, \quad \gamma = 0 = \bbmatrix{0&0\\0&0},
\]
then there created solutions $y=\bbmatrix{0&1}\gamma_*(x,t)\bbmatrix{1\\0}$ of the evolutionary Non Linear Schr\" odinger (ENLS) equation
\begin{equation} \label{eq:ENLS}
i y_t + y_{xx} + 2 |y|^2 y = 0.
\end{equation}
The proof of this fact can be found in \cite{bib:ENLS}.

Finally, solutions of the following canonical PDE
\begin{equation} \label{eq:CanSysEvol}
 \dfrac{\partial^2}{\partial t^2} \beta = \dfrac{\partial}{\partial x} [ - \dfrac{1}{2} (\beta')^2 - \dfrac{1}{4}\beta''' + \dfrac{(\dfrac{\partial}{\partial t} \beta)^2+ \dfrac{1}{4}(\beta'')^2}{\beta'}]
\end{equation}
are obtained if we choose
\[ \sigma_1 = \bbmatrix{0 & i \\ -i &  0},
   \sigma_2 = I = \bbmatrix{1 & 0 \\ 0 & 1},
   \gamma  = 0 = \bbmatrix{0 & 0 \\ 0 & 0}.
\]
The function $\beta(x,t) = \dfrac{\partial}{\partial x} \tau(x,t) = \dfrac{\partial}{\partial x} \det (\mathbb X^{-1}(0,0)\mathbb X(x,t))$
satisfies \eqref{eq:CanSysEvol}.

A reasonable question arises: are there other PDEs whose solutions are obtained using this scheme (i.e. solving vessel equations). To solve this question for the two-dimensions outer space ($\dim\mathcal E=2$), we define an obvious equivalence on vessel parameters. They are called equivalent if they produce the same first moment
$H_0(x,t)=C(x,t)\mathbb X^{-1}(x,t)B(x,t)$, which turns to be a solution of the same PDE. 
We show that there are just two major generalizations for the case $\mathcal E=\mathbb C^2$. One equation (for complex valued $\beta(x,t)$) is called a generalized KdV equation \eqref{eq:KdVGen}
\[ 4i \gamma_{12} (\beta)_t = 4(\gamma_{11}^2+\gamma_{12}\gamma_{21} )(\beta)_x -6 [(\beta)_x]^2-(\beta)_{xxx}, \]
and another equation is called a generalized ENLS equation \eqref{eq:NLSgen}
\[ 2a(\beta)_t = i \beta_{xx} -2i\gamma_{11} \beta_x + i2 (2a |\beta|^2 - \gamma_{21} \beta+\gamma_{12} h_{21}) (\gamma_{12} +2a\beta).
\]
Here, the parameters $a, \gamma_{ij}$ are constant numbers.

Future research will include other evolution equations for the variable $t$, creating corresponding hierarchies. A similar result already exists for the classical KdV equation \eqref{eq:KdV} in \cite{bib:KdVHierarchy}. Choosing $\mathcal E=\mathbb C^3$, we will obtain a family of the three dimensional outer space equations. For example, one of such vessels was used by the author to solve a Boussinesq equation in \cite{bib:BoussVessels}. Studying these classifications for all dimensions of $\mathcal E$,
we will be able to address the question which general PDEs are integrable.

Finally, we would like to mention some works of the author, originating the theory of vessels \cite{bib:SLVessels, bib:CanSys, bib:AnalPotKdV, bib:ENLS}  and joint works with collaborators \cite{bib:amv, bib:SchurIEOT, bib:MVZeroPole}.

\subsection{Vessels}
Let as start by presenting a more formal setting of vessels used to study solutions of completely integrable PDEs.
\begin{defn}
\label{def:VesPar}
A collection of $E\times E$ matrices $\sigma_1$, $\sigma_2$ and $\gamma$, where $\sigma_1$ is
invertible, are called vessel parameters.
\end{defn}
Vessel parameters are used to solve operator valued differential equations in the definition of a vessel.
\begin{defn} 
A \textbf{regular vessel} associated to the vessel parameters is a collection of bounded operators and spaces
\begin{equation} \label{eq:DefV}
\mathfrak{V} = (A_\zeta, C(x,t), \mathbb X(x,t), B(x,t), A; \sigma_1, 
\sigma_2, \gamma, \gamma_*(x,t);
\mathcal K,\mathcal E,\Omega),
\end{equation}
where $\mathcal K$ is a Krein space, $\mathcal E$ is a finite dimensional space of degree $E$ and $\sigma_1, \sigma_2, \gamma$ are vessel parameters. The following equations are assumed to hold:
\begin{subequations} \begin{align}
\label{eq:DBx} 0  =  \frac{d}{dx} (B(x,t)) \sigma_1 + A B(x,t) \sigma_2 + B(x,t) \gamma, \\
\label{eq:DBt}  \dfrac{\partial}{\partial t} B(x,t) = i A \dfrac{\partial}{\partial x} B(x,t).
\end{align}  \end{subequations}
\begin{subequations} \begin{align}
\label{eq:DCx} 0  = \sigma_1 \frac{d}{dx} (C(x,t)) + \sigma_2 C(x,t) A_\zeta - \gamma C(x,t), \\
\label{eq:DCt}  \dfrac{\partial}{\partial t} C(x,t) = -i \dfrac{\partial}{\partial x} C(x,t) A_\zeta.
\end{align}  \end{subequations}
\begin{subequations} \begin{align}
\label{eq:DXx} \frac{d}{dx} \mathbb X(x,t)  =  B(x,t) \sigma_2 C(x,t), \\
\label{eq:DXt} \dfrac{\partial}{\partial t} \mathbb X(x,t) =  i A B(x,t) \sigma_2 C(x,t) - i B(x,t)\sigma_2 C(x,t) A_\zeta + iB(x,t) \gamma C(x,t).
\end{align}  \end{subequations}
\begin{align}
\label{eq:Lyapunov} A \mathbb X(x,t) + \mathbb X(x,t) A_\zeta + B(x,t) \sigma_1 C(x,t) = 0, \\
\label{eq:Linkage}
\gamma_*(x,t)  = \sigma_2 C(x,t) \mathbb X^{-1}(x,t) B(x,t) \sigma_1
- \sigma_1 C(x,t) \mathbb X^{-1}(x,t) B(x,t) \sigma_2 .
\end{align}
where the operator $\mathbb X(x)$ is invertible in a region $\Omega$. Equation \eqref{eq:Lyapunov} is called \textbf{Lyapunov equation}, and \eqref{eq:Linkage}  is called the \textbf{linkage condition}.
\end{defn}
By definition, the transfer function of this vessel is
\begin{equation}\label{eq:DefS} 
S(\lambda,x,t) = I - C(x,t) \mathbb X^{-1}(x,t) (\lambda I - A)^{-1} B(x,t) \sigma_1,
\end{equation}
Poles and singularities of $S$ with respect to $\lambda$ are determined by $A$ only. In other words, it is a function possessing a Painleve property: its singularities are $x,t$ independent. Using Schur complements, the inverse function is given by
\begin{equation}\label{eq:DefS-1} 
S^{-1}(\lambda,x,t) = I + C(x,t) (\lambda I + A_\zeta)^{-1} \mathbb X^{-1}(x,t) B(x,t) \sigma_1.
\end{equation}
The vessel creates a completely integrable overdetermined input/state/output ($\mathrm u, \mathrm x, \mathrm y$) system as follows
\[ \Sigma_{2,3}: 
\left\{ \begin{array}{lll}
    \lambda \mathrm x_\lambda(x,t) = A  ~\mathrm x_\lambda(x,t) +  B(x,t) \sigma_1 ~\mathrm u_\lambda(x,t), \\[5pt]
    \dfrac{\partial}{\partial x} \mathrm x_\lambda(x,t) =  B(x,t) m_1 \mathrm u_\lambda(x,t), \\[5pt]
    \dfrac{\partial}{\partial t} \mathrm x_\lambda(x,t) = i A \dfrac{\partial}{\partial x} \mathrm x_\lambda(x,t) + i B(x,t) \sigma_1 \dfrac{\partial}{\partial x}\mathrm u_\lambda(x,t), \\[5pt]
    \mathrm y_\lambda(x,t) = \mathrm u_\lambda(x,t) - C(x,t)\mathbb X^{-1}(x,t)~ \mathrm x_\lambda(x,t).
  \end{array} \right.
\]
The system $\Sigma_{2,3}$ maps solutions $\mathrm u_\lambda(x,t)$ of
\begin{equation} \label{eq:InCCKdVtype}
\left\{ \begin{array}{ll}
\sigma_1 \dfrac{\partial}{\partial x} \mathrm u_\lambda(x,t) = \lambda m_1 \mathrm u_\lambda(x,t) + m_0 \mathrm u_\lambda(x,t), \\
\dfrac{\partial}{\partial t} \mathrm u_\lambda(x,t)  = -i\lambda \dfrac{\partial}{\partial x} \mathrm u_\lambda(x,t),
\end{array} \right. 
\end{equation}
to solutions $\mathrm y_\lambda(x,t) = S(\lambda,x,t)  \mathrm u_\lambda(x,t) $ of
\begin{equation} \label{eq:OutCCKdVtype}
\left\{ \begin{array}{ll}
\sigma_1 \dfrac{\partial}{\partial x} \mathrm y_\lambda(x,t) = \lambda m_1 \mathrm y_\lambda(x,t) + \gamma_*(x,t) \mathrm y_\lambda(x,t), \\
\sigma_1 \dfrac{\partial}{\partial t} \mathrm y_\lambda(x,t) = - i\lambda \dfrac{\partial}{\partial x} \mathrm y_\lambda(x,t) - i \sigma_1 \dfrac{\partial}{\partial x} [H_0(x,t)] \sigma_1 \mathrm y_\lambda(x,t),
\end{array} \right. 
\end{equation}
where $H_0(x,t) = C(x,t)\mathbb X^{-1}(x,t) B(x,t)$. In other words, the multiplication by $S(\lambda,x,t)$ serves as a B\" acklund transformations between equations of the wave type: if we carefully examine the system \eqref{eq:InCCKdVtype} it is easy to see that the first equation is a linear differential equation with constant coefficients and the second one is a wave equation. This fact can be found in \cite{bib:AnalPotKdV, bib:UnboundedVessels, bib:KdVVessels}, or alternatively an interested reader can prove the following, purely computational theorem.
\begin{thm}[Vessel=B\" acklund transformation] Assume that $u_\lambda(x,t)$ satisfies \eqref{eq:InCCKdVtype}. Let $\mathfrak V$ \eqref{eq:DefV} be a vessel and
$S(\lambda,x,t)$ \eqref{eq:DefS} its transfer function. Then $\mathrm y_\lambda(x,t) = S(\lambda,x,t)  \mathrm u_\lambda(x,t)$ satisfies \eqref{eq:OutCCKdVtype}.
\end{thm}
Notice that the existence of the function $H_0(x,t)$ requires invertability of the operator $\mathbb X(x,t)$. It is shown in \cite{bib:AnalPotKdV} in a more general setting that the operator $\mathbb X^{-1}(0,0)\mathbb X(x,t)$ possesses determinant. As a result, its invertability, and hence the invertability of the operator $\mathbb X(x,t)$ is equivalent to the determinant to be non-zero.
It is time to reveal the name of its determinant.
\begin{defn}\cite{bib:SLVessels} \label{def:Tau} The tau function $\tau(x,t)$ of the vessel $\mathfrak V$ is
\begin{equation} \label{eq:Tau} \tau = \det (\mathbb X^{-1}(0,0) \mathbb X(x,t)).
\end{equation}
\end{defn}
Developing the transfer function \eqref{eq:DefS}
\[
S(\lambda,x) = I - C\mathbb X^{-1}(\lambda I - A)^{-1}B\sigma_1 = I - \sum_{n=0}^\infty \dfrac{ H_n(x)}{\lambda^{n+1}}\sigma_1,
\]
into Taylor series, we arrive to the notion of moment: the \textbf{$n$-th moments} of a regular vessel $\mathfrak{V}$ is
$H_n(x) = C(x) \mathbb X^{-1}(x) A^n B(x)$. 
Moments play a key role in this research and we present some of their basic properties.
\begin{thm}[\cite{bib:KdVHierarchy, bib:AnalPotKdV}] The following relations between the moments of a KdV vessel $\mathfrak{V}$ hold
\begin{eqnarray}
\label{eq:DHnx} 	\sigma_1^{-1}\sigma_2H_{n+1} - H_{n+1} \sigma_2\sigma_1^{-1} &=& (H_n)'_x - \sigma_1^{-1} \gamma_* H_n + H_n \gamma\sigma_1^{-1}, \\
\label{eq:DHnt} 	(H_n)'_t& =& i (H_{n+1})'_x + i (H_0)'_x \sigma_1 H_n.
\end{eqnarray}
\end{thm}
Finally, the evolution of the generalized potential is as follows:
\begin{equation} \label{eq:Dgamma*tKdV}
(\gamma_*)'_t = - i \gamma_* (H_0)'_x\sigma_1 + i \sigma_1 (H_0)''_{xx} \sigma_1 +i \sigma_1 (H_0)'_x \gamma_*.
\end{equation} 
In order to see it, consider the mixed partial derivatives applied to the output of the system \eqref{eq:OutCCKdVtype}.

\subsection{Transformations of equivalency}
We discuss two types of transformations. One of them changes the output parameters and is called
external transformation. The other type changes the inner parameters and is called an inner transformation.
\begin{defn}
Two vessels producing solutions of the same PDE are called \textbf{equivalent}.
\end{defn}
Using the external transformation of equivalency, we all completely characterize KdV vessels with two dimensional outer space. These transformations may be further used for the higher dimensional outer space cases, which in turn might have additional transformations of equivalency.
\paragraph{External transformation I.}
Suppose that the collection \eqref{eq:DefV}
\[ \mathfrak{V} = (A_\zeta, C(x,t), \mathbb X(x,t), B(x,t), A; \sigma_1, 
\sigma_2, \gamma, \gamma_*(x,t);
\mathcal K,\mathcal E,\Omega),
\]
is a vessel. Define a new vessel
\begin{equation} \label{eq:DeftildeV}
\widetilde{\mathfrak{V}} = (A_\zeta, \widetilde C(x,t), \mathbb X(x,t), \widetilde B(x,t), A; \widetilde \sigma_1, 
\widetilde \sigma_2, \widetilde\gamma, \widetilde\gamma_*(x,t);
\mathcal K,\mathcal E,\Omega),
\end{equation}
where for two invertible $(\dim \mathcal E)^2$ matrices $V, U$
\begin{equation} \label{eq:OuterChange}
\begin{array}{ll}
\widetilde B(x,t) = B(x,t) U^{-1},  \quad \widetilde C(x,t) = V^{-1} C(x,t)\\
\widetilde \sigma_1 = U \sigma_1 V, \quad \widetilde \sigma_2 = U \sigma_2 V, \quad
\widetilde \gamma = U \gamma V, \quad \widetilde \gamma_*(x,t) = U \gamma_*(x,t) V.
\end{array} \end{equation}
Then immediate and easy computations show that the collection $\widetilde{\mathfrak{V}}$ is again a vessel. Such 
a vessel has exactly the same inner parameters: $\mathcal H, A_\zeta, \mathbb X(x,t), A$.

Taking the inverses of $V, U$ and using them for another external transformation, we obtain that the relation is actually an equivalence.
\begin{defn} \label{def:ExtEqI}Vessels $\mathfrak{V}$ and $\widetilde{\mathfrak{V}}$, defined in \eqref{eq:DefV} and \eqref{eq:DeftildeV}
respectively, are called \textbf{externally equivalent of the first kind}, if there exists a constant  $(\dim \mathcal E)^2$ invertible matrices $U, V$ such that
\eqref{eq:OuterChange} holds.
\end{defn}
How does the vessel ingredients change under equivalence relation? This is a very easy question to answer.
\begin{lemma}\label{lemma:ExtEqI} Assume that the vessels $\mathfrak{V}$ and $\widetilde{\mathfrak{V}}$, defined in \eqref{eq:DefV} and \eqref{eq:DeftildeV} respectively are externally equivalent. Then $\widetilde H_0(x,t) = H_0(x,t)$.
\end{lemma}
\begin{cor} Assume that the vessels $\mathfrak{V}$ and $\widetilde{\mathfrak{V}}$ are externally equivalent. Then
\[ \gamma_*(x,t)=\gamma_*(x,t), \quad \widetilde \tau(x)=\tau(x).
\]
\end{cor}
\Proof Immediate from Lemma \ref{lemma:ExtEqI} using definitions.  \qed

\paragraph{External transformation II.}
A more sophisticated transformation is presented next and will be used to classify PDEs, obtained from KdV Vessels.
\begin{thm}
Assume that vessels $\mathfrak V$ \eqref{eq:DefV} and $\widetilde{\mathfrak V}$ \eqref{eq:DeftildeV} have the following relation
between parameters
\begin{equation} \label{eq:S1Perturb}
\widetilde \sigma_2 = \sigma_2 + k_2 \sigma_1, \quad \widetilde \gamma = \gamma + k \sigma_1,
\end{equation}
then $\widetilde H_0(x,t) = H_0(x,t)$.
\end{thm}
\Proof  Equation \eqref{eq:DBx} with the new parameters $\widetilde \sigma_2, \widetilde \gamma$ becomes
\[ \dfrac{\partial}{\partial x} B \sigma_1 + A B \sigma_2 + k_2 A B \sigma_1 + B \gamma + k B \gamma = 0
\]
from which it is easy to deduce from variation of coefficients that
\[ \widetilde B(x,t) = V(x,t) B(x,t), \quad V(x,t)=\exp((A k_2 + k)(x+iAt)),
\]
where $B(x,t)$ is the solution of \eqref{eq:DBx} for the old parameters $\sigma_2, \gamma$. Similarly, one can
obtain that
\[ \widetilde C(x,t) = \widetilde C(x,t) U(x,t),
\]
where $U(x,t) = \exp((A_\zeta k_2 - k)(x-iAt))$. Finally, plugging
\[ \widetilde{\mathbb X}(x,t) = V(x,t) \mathbb X(x,t) U(x,t)
\]
we will obtain that new \eqref{eq:DXx}, and \eqref{eq:Lyapunov} hold.  As a result,
\[ \begin{array}{lll}
\widetilde H_0(x,t) & = \widetilde C(x,t) \widetilde{\mathbb X}^{-1} (x,t) \widetilde B(x,t) \\
& = C(x,t) U(x,t) (V(x,t) \mathbb X (x,t) U(x,t))^{-1} V(x,t) \widetilde B(x,t) \\
& = H_0(x,t). \qed
\end{array} \]
\begin{defn} Vessels $\mathfrak{V}$ and $\widetilde{\mathfrak{V}}$ are called \textbf{externally equivalent of the second kind}
if \eqref{eq:S1Perturb} holds.
\end{defn}

\paragraph{Internal transformation.}
It turns out that in order to produce solutions of PDEs one can focus on the simplest forms of the operators $A$, $A_\zeta$.
Suppose that $A=V^{-1} J V$ and $A_\zeta=U^{-1} J_\zeta U$, where $V, U$ are invertible operators and $J, J_\zeta$ are Jordan-block matrices, corresponding
to the eigenvalues. 
\begin{thm}\label{thm:V'} Assume that the collection
\[ \mathfrak{V} = (A_\zeta, C(x,t), \mathbb X(x,t), B(x,t), A; \sigma_1, 
\sigma_2, \gamma, \gamma_*(x,t);
\mathcal K,\mathcal E;\Omega)
\]
is a vessel. Then the collection
\[ \mathfrak{V}' = (A'_\zeta, C'(x,t), \mathbb X'(x,t), B'(x,t), A'; \sigma_1, 
\sigma_2, \gamma, \gamma_*(x,t);
\mathcal K,\mathcal E;\Omega),
\]
where
\begin{multline} \label{eq:V'Operators}
A'_\zeta = U^{-1} A_\zeta U, \quad C'(x,t) = C(x,t) U^{-1},  \quad \mathbb X'(x,t) = V \mathbb X(x,t) U^{-1}, \\
B'(x,t)=VB(x,t),  \quad A' = V A V^{-1}
\end{multline}
is also a vessel possessing the same moments and the same $\gamma_*(x,t)$.
\end{thm}
\Rem Notice that the operators $A, A_\zeta$ become Jordan-block matrices in this case.

\Proof Notice first that the equations \eqref{eq:DBx} and \eqref{eq:DBt} can be multiplied by $V$ on the left to produce the corresponding vessel equations for
$\mathfrak{V}'$:
\[ \begin{array}{llllll}
\frac{d}{dx} (B'(x,t)) + A' B'(x,t) \sigma_2 + B'(x,t) \gamma = \frac{d}{dx} (V B(x,t)) + V A V^{-1} V B(x,t) \sigma_2 + V B(x,t) \gamma = 0, \\
\dfrac{\partial}{\partial t} B'(x,t) = \dfrac{\partial}{\partial t} VB(x,t) = i VAV^{-1} \dfrac{\partial}{\partial x} VB(x,t) = i A' \dfrac{\partial}{\partial x} B'(x,t).
\end{array} \]
In the same way one obtains equations \eqref{eq:DCx}, \eqref{eq:DCt} by multiplying the same equations for $\mathfrak{V}$ from the right on $U^{-1}$. Finally,
equations \eqref{eq:DXx}, \eqref{eq:DXt} and \eqref{eq:Lyapunov} for $\mathfrak{V}'$ are obtained from the same equations for $\mathfrak{V}$ by
multiplying them on $V$ from the left and on $U^{-1}$ from the right.
\qed

Once the simplest form of the operators $A, A_\zeta$ is established, there is a question arises regarding the form of the function-matrices $B, C$ and $\mathbb X$.
The only freedom that is remained is in further multiplication by matrices $V, U$ that commute with $A$ and $A_\zeta$ respectively. This is done in order not to change the
established operators $A, A_\zeta$ in their simplest form.
\begin{thm} \label{thm:AVVA}Assume that $A V = V A$, $U A_\zeta = A_\zeta U$, then if the collection
\[ \mathfrak{V} = (A_\zeta, C(x,t), \mathbb X(x,t), B(x,t), A; \sigma_1, 
\sigma_2, \gamma, \gamma_*(x,t);
\mathcal K,\mathcal E;\Omega)
\]
is a vessel, so is the collection
\[ \mathfrak{V}' = (A_\zeta, C'(x,t), \mathbb X'(x,t), B'(x,t), A; \sigma_1, 
\sigma_2, \gamma, \gamma_*(x,t);
\mathcal K,\mathcal E;\Omega),
\]
where
\[ C'(x,t) = C(x,t) U^{-1},  \quad \mathbb X'(x,t) = V \mathbb X(x,t) U^{-1}, \quad B'(x,t) = V B(x,t).
\]
\end{thm}
\Proof similar to the proof of Theorem \ref{thm:V'}. \qed
\begin{defn} Vessels $\mathfrak{V}$ and $\mathfrak{V}'$ are called \textbf{internally equivalent}
if \eqref{eq:V'Operators} holds.
\end{defn}

\begin{mthm} Assume that 
\[ \mathfrak{V} = (A_\zeta, C(x,t), \mathbb X(x,t), B(x,t), A; \sigma_1, 
\sigma_2, \gamma, \gamma_*(x,t);
\mathcal K,\mathcal E;\Omega)
\]
is a vessel. Then $\mathfrak V$ is equivalent to the following vessel
\[ \mathfrak{V} = (J_\zeta, C(x,t), \mathbb X(x,t), B(x,t), J; I, j_2, \gamma, \gamma_*(x,t);
\mathcal K,\mathcal E;\Omega)
\]
where $J,J_\zeta$ are in simplest via similarity form and $j_2$ is a Jordan block matrix.
\end{mthm}
\Proof from the external transformation of the first kind, we will obtain that $\sigma_1=I$ and applying it again, that $\sigma_2$ is in a Jordan-Block form. Then using the internal transformation, the operators $A,A_\zeta$ are brought to the simplest via similarity form.
\qed
\section{Recurrence relations on moment entries}
let us further investigate formula \eqref{eq:DHnx}
\[ \sigma_1^{-1}\sigma_2H_{n+1} - H_{n+1} \sigma_2\sigma_1^{-1} = (H_n)'_x - \sigma_1^{-1} \gamma_* H_n + H_n \gamma\sigma_1^{-1}. \]
Taking matrix $K$ satisfying $[K,\sigma_1^{-1}\sigma_2] = K \sigma_1^{-1}\sigma_2 - \sigma_1^{-1}\sigma_2 K$, and multiplying \eqref{eq:DHnx} by $\sigma_1$ 
we will obtain that
\[ 0 = \TR(K [\sigma_1^{-1}\sigma_2, H_{n+1}\sigma_1] ) = \TR(K (H_n)'_x \sigma_1 - K\sigma_1^{-1} \gamma_* H_n\sigma_1 + K H_n \gamma)
\]
or
\[ \TR(K (H_n)'_x \sigma_1) = \TR((\sigma_1 K\sigma_1^{-1} \gamma_* -\gamma K) H_n).
\]
On the other hand, if $[K,\sigma_1^{-1}\sigma_2]\neq 0$, we will obtain a relation on the entries of $H_{n+1}$ in terms of the previous moment $H_n$ and $\gamma_*$.
Thus the following Lemma holds.
\begin{lemma} If moments satisfy equation \eqref{eq:DHnx} then for any $E\times E$ matrix $K$ it holds
\begin{eqnarray}
\label{eq:Kcomm} K \sigma_1^{-1}\sigma_2 = \sigma_1^{-1}\sigma_2 K & \Rightarrow & \TR(K (H_n)'_x \sigma_1) = \TR((\sigma_1 K\sigma_1^{-1} \gamma_* -\gamma K) H_n) \\
\label{eq:KNotcomm} K \sigma_1^{-1}\sigma_2 \neq \sigma_1^{-1}\sigma_2  K & \Rightarrow & \TR((K \sigma_1^{-1}\sigma_2 -\sigma_1^{-1}\sigma_2 K) H_{n+1}\sigma_1 ) =  \\ 
\nonumber && \hspace{1cm} = \TR(K (H_n)'_x \sigma_1 - K\sigma_1^{-1} \gamma_* H_n\sigma_1 + K H_n \gamma).
\end{eqnarray}
\end{lemma}
It follows from this lemma that if there is a matrix $K = [\sigma_1^{-1}\sigma_2,L]$ for some $L$ and $K$ additionally commutes with $\sigma_1^{-1}\sigma_2$
then
\[ \begin{array}{lll}
\TR(K (H_{n+1})'_x \sigma_1) & =  \text{ by \eqref{eq:Kcomm} } = \TR((\sigma_1 K\sigma_1^{-1} \gamma_* -\gamma K) H_{n+1}) \\
	& = \text{ by $\dfrac{\partial}{\partial x}$\eqref{eq:KNotcomm} } = - \dfrac{\partial}{\partial x} \TR(L (H_n)'_x \sigma_1 - L\sigma_1^{-1} \gamma_* H_n\sigma_1 + L H_n \gamma).
\end{array} \]
Plugging here the formula \eqref{eq:Linkage} for $\gamma_*$ we can obtain that
\[ \TR([K,\sigma_1^{-1}\gamma]H_{n+1}\sigma_1) = \dfrac{\partial}{\partial x} \TR(L (H_n)'_x \sigma_1 - L\sigma_1^{-1} \gamma_* H_n\sigma_1 + L H_n \gamma)-
\TR(KH_0\sigma_1((H_n)'_x \sigma_1 - \sigma_1^{-1} \gamma_* H_n\sigma_1 + H_n \gamma))
\]
\section{Canonical parameters}
Using the external equivalency of the first kind we can always assume that $\sigma_1=I$. Indeed, choosing
$U=\sigma_1^{-1}$, $V=I$ in transformations \eqref{eq:OuterChange} we will obtain that. Of course, the parameters
$\sigma_2, \gamma$ will be changed too. Next step is to bring the parameter $\sigma_2$ to its Jordan-block form, preserving the fact that $\sigma_1=I$. Indeed assuming $\sigma_1=I$, $\sigma_2=U^{-1}JU$ we can use transformations \eqref{eq:OuterChange} with $V=U^{-1}$ achieving $\widetilde\sigma_1=U^{-1} I U=I$, $\widetilde\sigma_2 = J$. Finally using the external transformation of the second kind, we will obtain that any set of parameters is equivalent to the following:
\begin{defn} A set of $2\times 2$ vessel parameters $\sigma_1, \sigma_2, \gamma$ is called \textbf{canonical} if $\sigma_1=I$,
$\gamma=\bbmatrix{\gamma_{11}&\gamma_{12}\\\gamma_{21}&-\gamma_{11}}$ (i.e. $\TR(\gamma)=0$) and either
\begin{itemize}
\item $\sigma_2=\bbmatrix{a&0\\0&-a}, a\neq 0$, or
\item $\sigma_2=\bbmatrix{0&0\\1&0}$.
\end{itemize}
In the first case, the triple $\sigma_1=I$, $\bbmatrix{a&0\\0&-a}, a\neq 0$, and $\gamma=\bbmatrix{\gamma_{11}&\gamma_{12}\\\gamma_{21}&-\gamma_{11}}$ is called \textbf{generalized NLS parameters}. In the second case, the triple
$\sigma_1=I$, $\sigma_2=\bbmatrix{0&0\\1&0}$, and $\gamma=\bbmatrix{\gamma_{11}&\gamma_{12}\\\gamma_{21}&-\gamma_{11}}$ is called \textbf{generalized KdV vessel parameters}.
\end{defn}
\begin{mthm} Any set of $2\times 2$ vessel parameters is equivalent to canonical vessel parameters.
\end{mthm}
\Proof We have seen that using external equivalence of the first kind, we can assume that $\sigma_1=I$ and $\sigma_2=J$, its Jordan-block form. Moreover, using the external transformation of the second kind, we immediately obtain that $\TR(\gamma)=0$ and 
$\sigma_2$ has also trace zero, in the case of two distinct eigenvalues, or that the eigenvalue can be removed by the second external transformation, leaving only $1$ in the off-diagonal entry (one of them, which we choose the 2,1 entry without loss of generality). 
\qed

It follows from these considerations that canonical systems, presented in \cite{bib:CanSys} are equivalent to the generalized ENLS equation. For the canonical systems one uses parameters
\[ \sigma_1=\bbmatrix{0&i\\-i&0}, \quad \sigma_2=I,\quad \gamma=0.
\]
Consequently, using the external transformation of the first kind for $U=\sigma_1$, $V=I$ we will obtain that these parameters are equivalent to
\[ \sigma_1=I, \quad \sigma_2=\bbmatrix{0&i\\-i&0}, \quad \gamma=0.
\]
Finally, since $\sigma_2=V \bbmatrix{1&0\\0&-1} V^{-1}$ for $V=\dfrac{1}{2} \bbmatrix{i&-i\\1&1}$ we conclude that, again from the external transformation of the first kind for $U=V$ these vessel parameters are equivalent to the canonical generalized ENLS parameters with $\gamma=0$.

\subsection{generalized NLS parameters}
Let us develop the formula, corresponding to generalized NLS parameters.  
Assume that $H_0(x,t)=\bbmatrix{h_{11}&h_{12}\\h_{21}&h_{22}}$, then
\[ \begin{array}{lll}
\gamma_*(x,t) & = \gamma + \sigma_2 H_0(x,t) \sigma_1 - \sigma_1 H_0(x,t) \sigma_2 \\
& = \gamma + \bbmatrix{a&0\\0&-a} \bbmatrix{ h_{11} & h_{12} \\ h_{21} & h_{22}} - \bbmatrix{h_{11}&h_{12}\\h_{21}&h_{22}} \bbmatrix{a&0\\0&-a} \\
& = \gamma + \bbmatrix{0&(a-(-a))h_{12}\\(-a-a)h_{21} & 0} \\
& = \bbmatrix{\gamma_{11} & \gamma_{12} +2ah_{12} \\ \gamma_{21} -2a h_{21} & -\gamma_{11} }.
\end{array} \]
The evolutionary equation \eqref{eq:Dgamma*tKdV} is equivalent to the evolution of its $1,2$ entry $(b-a)h_{12}$  and we will obtain that
\[ \begin{array}{llllll}
\bbmatrix{1&0}(\gamma_*)_t \bbmatrix{0\\1} & = \bbmatrix{1&0}[-i \gamma_* (H_0)_x \sigma_1 + i \sigma_1 (H_0)_{xx} \sigma_1 +i \sigma_1 (H_0)_x \gamma_* ]\bbmatrix{0\\1}\\
2a (h_{12})_t & = -i \bbmatrix{\gamma_{11} &\gamma_{12} +2ah_{12}} \bbmatrix{h'_{12}\\h'_{22}} + i h''_{12} + i \bbmatrix{h'_{11}&h'_{12}} \bbmatrix{\gamma_{12} +2ah_{12}\\-\gamma_{11} } \\
		& = i h''_{12} -i\gamma_{11}  h_{12}'-i(\gamma_{12} +2ah_{12})h'_{22}+ih'_{11}(\gamma_{12} +2ah_{12})-i\gamma_{11} h'_{12}\\
		& = i h''_{12} -2i\gamma_{11} h_{12}' + i(h'_{11}-h'_{22})(\gamma_{12} +2ah_{12}).
\end{array} \]
In order to find $h'_{11}-h'_{22}$ in terms of $h_{12}$ we examine the formula \eqref{eq:DHnx}:
\[ \begin{array}{ll}
\sigma_2 H_1 - H_1 \sigma_2 & = H_0'-\gamma_* H_0+H_0\gamma \\
& = H_0'-\bbmatrix{\gamma_{11} & \gamma_{12} +2ah_{12} \\ \gamma_{21} -2a h_{21} & -\gamma_{11} } 
\bbmatrix{h_{11}&h_{12}\\h_{21}&h_{22}} + \bbmatrix{h_{11}&h_{12}\\h_{21}&h_{22}} \bbmatrix{\gamma_{11}&\gamma_{12}\\\gamma_{21}&-\gamma_{11}}
\end{array} \]
and take $1,1$ and $2,2$ entries on both sides:
\[ \begin{array}{ll}
0 = h'_{11} - 2a h_{12}h_{21} - \gamma_{12} h_{21}+\gamma_{21} h_{12},  \\
0 = h'_{22} + 2a h_{12}h_{21} - \gamma_{21} h_{12}+\gamma_{12} h_{21}.
\end{array} \] 
Thus
\[ h'_{11}-h'_{22} = 2 (2a h_{12}h_{21} - \gamma_{21} h_{12}+\gamma_{12} h_{21})
\]
and plugging it back into the evolutionary equation, we will obtain the following formula:
\begin{equation}\label{eq:NLSgen}
2a(h_{12})_t = i h''_{12} -2i\gamma_{11} h_{12}' + i2 (2a h_{12}h_{21} - \gamma_{21} h_{12}+\gamma_{12} h_{21}) (\gamma_{12} +2ah_{12})
\end{equation}
This is the standard NLS equation in case $\gamma_{11}=0$, $a=-1/2$, and $h_{21}=h_{12}^*$.

\subsection{Generalized KdV parameters}
Consider next the second canonical case. Let as denote  $\gamma=\bbmatrix{\gamma_{11}&\gamma_{12}\\\gamma_{21}&-\gamma_{11}}$ and  $H_0(x,t)=\bbmatrix{h_{11}&h_{12}\\h_{21}&h_{22}}$. Then 
\[ \begin{array}{lll}
\gamma_*(x,t) & = \gamma + \sigma_2 H_0(x,t) \sigma_1 - \sigma_1 H_0(x,t) \sigma_2 \\
& = \gamma + \bbmatrix{0&0\\1&0} \bbmatrix{h_{11}&h_{12}\\h_{21}&h_{22}}  -  \bbmatrix{h_{11}&h_{12}\\h_{21}&h_{22}} \bbmatrix{0&0\\1&0} \\
& = \bbmatrix{\gamma_{11}-h_{12} &\gamma_{12}\\\gamma_{21}+h_{11}-h_{22} &-\gamma_{11}+h_{12}}.
\end{array} \]
Plugging $K=\sigma_1$ and then $ \sigma_2$ into \eqref{eq:Kcomm} we will obtain that 
\[ (h_{11})_x+(h_{22})_x = 0.
\]
and
\[  (h_{12})_x = \gamma_{12}(h_{22}-h_{11}) + (\gamma_{11}-\gamma_{22})h_{12}-h_{12}^2
\]
From where it follows that under condition $\gamma_{12}\neq 0$
\[ h_{22}=-h_{11}, \quad h_{11} = - \dfrac{-2 \gamma_{11} h_{12} + h_{12}^2 + (h_{12})_x}{2\gamma_{12}}.
\]
Then taking the 1,1 entry of the equality \eqref{eq:Dgamma*tKdV} leads to
\begin{equation}\label{eq:KdVPre} \begin{array}{llll}
-(h_{12})'_t & = - \bbmatrix{1&0}(i \gamma_* (H_0)'_x\sigma_1 + i \sigma_1 (H_0)''_{xx} \sigma_1 +i \sigma_1 (H_0)'_x \gamma_*)\bbmatrix{1\\0} \\
& = -\frac{i}{2\gamma_{12}} (2 (h_{12})_x \left(-2 \gamma_{11} h_{12} +h_{12}^2-\gamma_{12} \gamma_{21}\right)+2 (h_{12})_{xx} (h_{12}-\gamma_{11}) + \\
& \hspace{2cm} + 4 [(h_{12})_x]^2 + (h_{12})_{xxx} + 2 \gamma_{12}^2 (h_{21})_x ).
\end{array} \end{equation}
There is an unknown at this point parameter $(h_{21})_x$, which can be found from the equations \eqref{eq:Kcomm}, \eqref{eq:KNotcomm} as follows. Plugging $K=\bbmatrix{1&0\\0&0}$ and then $\gamma$ into \eqref{eq:KNotcomm}, we will obtain formulas for the entries of $H_1(x,t)$. If we denote the elements of $H_1(x,t) = \bbmatrix{k_{11}&k_{12}\\k_{21}&k_{22}}$, they are
\[ \begin{array}{lll}
 k_{22} = \dfrac{1}{2\gamma_{12}^2} \big( 2\gamma_{12}^2 k_{11}-4\gamma_{11}\gamma_{12}k_{12} -2\gamma_{11} h_{12}^3 +h_{12}^4 -(h_{12})_x (4\gamma_{11}^2 + 2\gamma_{12}\gamma_{21}) + [(h_{12})_x]^2 + \\
 \hspace{1.5cm} +2 h_{12}^2(-\gamma_{12}\gamma_{21} +(h_{12})_x) +2 h_{12} (\gamma_{12} h_{21} +\gamma_{11} (h_{12})_x) -
 2 \gamma_{12}^2 (h_{21})_x +2\gamma_{11} (h_{12})_{xx} \big),\\
k_{12} = \dfrac{1}{2\gamma_{12}} \big( (h_{12})_{xx} +h_{12}(3(h_{12})_x -2\gamma_{12}\gamma_{21}) -2\gamma_{11} (h_{12})_x + 2\gamma_{12}^2 h_{21} + h_{12}^3 -2\gamma_{11}h_{12}^2  \big).
\end{array} \]
Consider \eqref{eq:Kcomm} for $K=\sigma_2$, substituting there the expressions for $k_{22}, k_{12}$ just found. After simple cancellations we will arrive to
\[ (h_{21})_x = -\dfrac{1}{2 \gamma_{12}^2 } \big( 2 (h_{12})_x \left(-2 \gamma_{11} h_{12}+h_{12}^2-\gamma_{12} \gamma_{21} \right) + 2 (h_{12})_{xx} (h_{12}-\gamma_{11})+4 [(h_{12})_x]^2+(h_{12})_{xxx} \big).
\]
Finally when we plug this found expression for $(h_{21})_x$ into \eqref{eq:KdVPre}, we will obtain
\begin{equation}\label{eq:KdVGen}
4i \gamma_{12} (h_{12})_t = 4(\gamma_{11}^2+\gamma_{12}\gamma_{21} )(h_{12})_x -6 [(h_{12})_x]^2-(h_{12})_{xxx},
\end{equation}
which is the famous KdV equation for $\gamma_{12}=i$ and $ \gamma_{11}=\gamma_{21}=0$. More precisely, to get the KdV equation one will need to differentiate the PDE with respect to $x$ and use the substitution $f=(h_{12})_x$.

\section{Soliton solutions}
Soliton solutions are obtained if we use a finite dimensional inner space $\mathcal K$. For the two canonical cases, we present a one dimensional solitons for the most general case, leaving to the reader to fill in the details, to explore some particular cases and higher dimensional solitons.
\paragraph{Generalized ENLS soliton.} Notice first that the differential equations \eqref{eq:DBx}, \eqref{eq:DCx} has constant coefficients and are solvable. To simplify computations,
we present a simpler case when $\gamma_{12}=\gamma_{21}=0$. In this case, we denote
\[ k = A a + \gamma_{11}, \quad k_\zeta = -A_\zeta a + \gamma_{11}
\]
and functions of the vessel as follows:
\[ \begin{array}{ll}
B(x,t) = \bbmatrix{\exp(-kx-iAkt) B1 & \exp(kx+iAkt) B2 }, \\
C(x,t) = \bbmatrix{\exp(k_\zeta x-iA_\zeta k_\zeta t) C1\\ \exp(-k_\zeta x+iA_\zeta k_\zeta t) C2},
\end{array} \]
where $B1, B2, C1, C2$ are arbitrary constants. Define next
\[ \mathbb X(x,t) = \dfrac{B(x,t) C(x,t)}{A+A_\zeta}.
\]
Then for the first moment $H_0(x,t)=C(x) \mathbb X^{-1}(x) B(x)=[h_{ij}]$ we will obtain that
\[ h_{12} = -\frac{C1 (A+A_\zeta) e^{2 (a A+\gamma_{11}) (x+i A t)}}{B1C1+C2 \exp(2 i (A+A_\zeta) (\gamma_{11} t+a (A t-A_\zeta t-i x))}
\]
and
\[ h_{21} = -\frac{B1 C2 (A+A_\zeta) e^{2 (a \text{Azeta}-\gamma_{11}) (x-i A_\zeta t)}}{B1C1+C2 \exp (2 i (A+A_\zeta) (\gamma_{11} t+a (A t-A_\zeta t-i x)))}.
\]
In order to get an equation in one function ($h_{12}$) we can choose symmetric parameters and realizations. For this it is enough to require
\[ A_\zeta=A^*, B1=1, C1=C2, a, i \gamma_{11}\in\mathbb R,
\]
and then $h_{21}=h_{12}^*$.

\paragraph{Generalized KdV soliton.} Suppose that $A, A_zeta$ are fixed numbers together with the $2\times 2$ matrix-parameter $\gamma=\bbmatrix{\gamma_{11}&\gamma_{12}\\\gamma_{21}&-\gamma_{11}}$. Define the following numbers
\[ k = \sqrt{A\gamma_{12}-\det\gamma}, \quad \quad k_\zeta=\sqrt{-A_\zeta \gamma_{12}-\det\gamma}
\]
and the functions (for arbitrary constants $B1, B2, C1, C2$)
\[ \begin{array}{lllll}
b_2(x,t) = \cosh(k x+ i A k t) B1 + \sinh(k x+ i A k t) B2, 
\quad b_1(x,t) = \dfrac{-(b_2(x,t))_x + \gamma_{11} b_2(x,t) }{\gamma_{12}} \\
B(x,t) = \bbmatrix{b_1(x,t)&b_2(x,t)} \\
c_1(x,t) = \cosh(k_\zeta x- i A_\zeta k t) C1 + \sinh(k_\zeta x- i A_\zeta k_\zeta t) C2, 
\quad c_2(x,t) = \dfrac{(c_1(x,t))_x - \gamma_{11} c_1(x,t) }{\gamma_{12}} \\
C(x,t) = \bbmatrix{c_1(x,t)\\c_2(x,t)}
\end{array} \]
and finally
\[ \mathbb X(x,t) = \dfrac{B(x,t) C(x,t)}{A+A_\zeta}.
\]
Then this collection defines a vessel, producing a one-dimensional soliton for the generalized KdV equation \eqref{eq:KdVGen}:
\[ h_{12} = \dfrac{\bbmatrix{1&0}C(x,t)B(x,t)\bbmatrix{0\\1}}{\mathbb X(x,t)}.
\]

\end{document}